\documentclass{article}
\usepackage{latexsym,amssymb,amsmath}
\usepackage{amsfonts}

\newtheorem{theorem}{Theorem}[section]
\newtheorem{lemma}[theorem]{Lemma}
\newtheorem{proposition}[theorem]{Proposition}
\newtheorem{corollary}[theorem]{Corollary}

\newtheorem{remark}[theorem]{Remark}
\newenvironment{proof}[1][Proof]{\textbf{#1.} }{\ \mathbb Rule{0.5em}{0.5em}}

\newcommand{\pa}{\partial}
\newcommand{\dom}{{\mathcal D}}

\begin{document}

\title{Some examples of absolute continuity of measures in stochastic
  fluid dynamics}
\author{B. Ferrario\\     Dipartimento di Matematica - Universit\`a di Pavia}
\date{}
\maketitle

\begin{abstract}
A non linear It\^o equation in a Hilbert space is studied by means of
Girsanov theorem. We consider a non linearity of polynomial growth
in suitable norms,
including that of quadratic type which
appears in the  Kuramoto--Sivashinsky equation and in the 
Navier--Stokes equation. We prove that Girsanov theorem holds for the
1-dimensional stochastic Kura\-moto--Sivashinsky equation  and for 
a modification of the 
2- and 3-dimensional stochastic Navier--Stokes equation; this
modification consists in substituting the 
Laplacian $-\Delta$ with $(-\Delta)^\alpha$, where $\alpha>\frac d2
+1$ ($d=2,3$).
In this way, we prove existence and uniqueness of
solutions for these stochastic equations.
Moreover, the asymptotic behaviour for  $t \to \infty$  is characterized.
\end{abstract}

\section{Introduction}\label{intro}
The study of non linear equations requires some skill
to deal successfully with the non linearity.
As fas as stochastic differential equations are concerned, a possible 
technique to study a non linear equation is the Girsanov transform.
Indeed, given a non linear stochastic It\^o
equation
\begin{equation}\label{nlin}
du(t)+\left[\ Lu(t) + F(u(t))\ \right] dt=Gdw(t), \quad t\in ]0,T];\; u(0)=x
\end{equation}
it is possible to analyze it 
as a perturbation of the linear equation 
\begin{equation}\label{lin}
dz(t)+Lz(t)\ dt=Gdw(t), \quad t\in ]0,T]; \; z(0)=x
\end{equation}
by means of Girsanov theorem.
It is well known that this theorem 
holds  if for instance Novikov condition 
$$
 \mathbb E \Big[\exp\big(\tfrac 12 \textstyle \int_0^T |G^{-1} F(z(t))|^2 dt\big) \Big]
   <\infty
$$
is satisfied.
We are interested in the case in which Novikov condition is not
fulfilled, but it is if
the non linear term $F$ is suitably truncated. By an approximating procedure 
we can apply Girsanov transform so to get that
equation \eqref{nlin} has a weak solution having the same regularity
as $z$ and  the law 
$\mathcal L_u$ of $u$ is absolutely continuous with respect to the law 
$\mathcal L_z$ of $z$ ($\mathcal L_u\prec \mathcal L_z$) 
and possibly the converse too, so to get the
equivalence of $\mathcal L_u$ and $\mathcal L_z$ 
($\mathcal L_u \sim \mathcal L_z$).
We recall that if $\mathcal L_u\prec \mathcal L_z$, 
uniqueness for equation \eqref{lin} 
implies uniqueness in law for equation \eqref{nlin}. Moreover, if
$\mathcal L_u \sim \mathcal L_z$, each property
holding P-a.s. for the
process $z$ must hold also for the process $u$ and vice versa.

Our analysis to verify if Girsanov transform can be used 
is quite standard. We formalize it here in order to apply
it in the next sections 
to some models in stochastic fluid dynamics, in which the equations
are set in a infinite dimensional Hilbert space and the non
linearity $F$ is of quadratic type; however, the result holds true for
$F$ of polynomial growth.

As to the structure of the paper, in Section \ref{leggi}
two abstract results are presented; in the first it is proved that
$\mathcal L_u\prec \mathcal L_z$
and in the second that    $\mathcal L_u \sim \mathcal L_z$.
Then, in the other two sections these results are applied  to
a stochastic Kuramoto--Sivashinsky
equation and to a modified stochastic Navier--Stokes equation, respectively.
For these non linear equations (which have a similar non linearity), 
we  obtain results of
existence and uniqueness of the solution;
further, uniqueness of the invariant measure is proved, so to
characterize the  asymptotic behaviour.

\section{Absolute continuity of laws}
             \label{leggi}
We are given a separable Hilbert space $H$, equipped 
 with a complete orthonormal
system $\{e_j\}_{j=1}^\infty$,
and a complete probability space $(\Omega, F,\{F_t\}_{t\ge 0},P)$.
We denote by $\mathbb E$ the expectation with respect to the measure $P$.

As far as equation \eqref{lin} is concerned,
we assume that $L$ and $G$ are linear operators in $H$ and $G$ is invertible.
The process $w$ is a cylindrical Wiener process in $H$,
defined on the probability space 
$(\Omega, F,\{F_t\}_{t\ge 0},P)$.
This means that, given a sequence $\{\beta_j\}_{j=1}^\infty$ 
of  i.i.d. one dimensional Wiener processes defined on
$(\Omega, F,\{F_t\}_{t\ge 0},P)$, we represent $w(t)=\sum_j
\beta_j(t) e_j$.
\\
Moreover, we assume that
there exists a unique strong solution $z$ (in the stochastic sense)
which is a Markov process such that 
\begin{equation}\label{z:2p}
 \mathbb E \|z\|^{2p}_{C([0,T];E)}<\infty
\end{equation}
for some $p>1$, 
where $E$ is a separable dense subset of $H$.
Actually, it would be enough $z$ to be a weak solution; but in our
applications in Sections \ref{1DeqKS} and \ref{navier-stokes}, 
$z$ will be a strong
solution and thus we assume it since now.

From now on, we denote by $z(t;x)$, or simply by $z(t)$,
 the solution of \eqref{nlin} evaluated
at time $t$ (thus $z(0;x)=x$ and, for $t>0$, $z(t;x)$ is a random variable) 
and by $z$ the solution process $\{z(t;x)\}_{0\le t \le T}$ on a time
interval $[0,T]$.


The main assumption on the non linear term is that the operator
$G^{-1}F:E\to H$ is  measurable and 
\begin{equation}\label{crescita}
|G^{-1}F(v)|_H\le c\left( 1+|v|^p_E\right)
 \qquad \forall v \in E,
\end{equation}
where $c$ is a suitable constant and $p>1$ is the same as in \eqref{z:2p}.
This implies that 
\begin{equation}\label{stima1}
 \int_0^T |G^{-1}F(z(t))|_H ^2 dt 
 \le 2 T c^2 \left( 1+\|z\|^{2p}_{C([0,T];E)}\right) 
\end{equation}
so
\begin{equation}\label{GF:2}
 \mathbb E  \int_0^T |G^{-1}F(z(t))|_H ^2 dt  
 \le
 2 T c^2 \left( 1+\mathbb E[\|z\|^{2p}_{C([0,T];E)}]\right)<\infty. 
\end{equation}
In particular
\begin{equation}\label{pathNov}
 P\big\{ \textstyle\int_0^T |G^{-1}F(z(t))|_H ^2 dt <\infty \big\}=1.
\end{equation}
This condition is necessary for Novikov condition 
\begin{equation}\label{nov}
 \mathbb E \Big[\exp\big(\tfrac 12 \textstyle\int_0^T |G^{-1} F(z(t))|_H^2 dt\big)\Big]
 <\infty
\end{equation}
to hold.
It is well known (see, e.g., \cite{dpz}
for stochastic PDE's in Hilbert spaces) that condition \eqref{nov} implies 
that 
\begin{equation}\label{densita}
\rho^T_{u/z}:=
 \mathbb E\left[\exp\Big(\textstyle\int_0^T \!_H\langle G^{-1}F(z(s)),dw(s)\rangle_H
 -\tfrac 12 \textstyle\int_0^T  |G^{-1} F(z(s))|_H^2ds \Big)
 \big|\sigma^T(z)\right] 
\end{equation}
is a probability  density.
Here $\sigma^T(z)$ denotes the $\sigma$-algebra generated by
$\{z(t)\}_{0\le t\le T}$. 
The stochastic integral in the exponent has to be understood as
$\sum_j \int_0^T \!_H\langle G^{-1}F(z(s)),e_j\rangle_H d\beta_j(s) $
and is well defined because of \eqref{GF:2}
(see \cite{dpz}, Chapter 4).

As soon as we know that $\mathbb E[\rho^T_{u/z}]=1$, we apply Girsanov
theorem to get that $\mathcal L_u\prec \mathcal L_z$. We remind it
here, for reader's convenience (see, e.g., \cite{dpz}, \cite{masl}, \cite{mr}).
Defined the probability measure $P^*$ on $(\Omega,F)$ 
by $dP^*=\rho^T_{u/z} dP$,  Girsanov theorem states that
$$
 w^*(t)=w(t)+\int_0^t G^{-1}F(z(s))ds
$$
is a cylindrical Wiener process on $(\Omega,F,\{F_t\}_{0\le t\le T}, P^*)$.
So, if $z$ solves equation \eqref{lin} with Wiener process $w$, then
$z$ solves equation \eqref{nlin} with Wiener process $w^*$, since
\begin{equation*}
\begin{split}
 z(t)&=x-\int_0^t Lz(s)ds+\int_0^t Gdw(s) 
  \\ &=x-\int_0^t Lz(s)ds-\int_0^t F(z(s))ds+\int_0^t Gdw^*(s).
\end{split}
\end{equation*}
Thus, $P\{u \in \Lambda\}=P^*\{z \in \Lambda \}$
for every Borel set $\Lambda\subset C([0,T];E)$.
Then  $P\{z \in \Lambda\}=0$ implies 
$P^*\{z \in \Lambda \}=0$ and so $P\{u \in \Lambda\}=0$, that is 
$\mathcal L_u \prec \mathcal L_z$.

Summing up, 
assuming that the solution $z$ to 
equation \eqref{lin} is such that $\mathbb E[\rho^T_{u/z}]=1$,
then 
equation \eqref{nlin} has a weak solution having the same regularity
as $z$ and  $\mathcal L_u \prec \mathcal L_z$;
moreover, uniqueness in law for $z$ implies uniqueness in law for $u$.
\\
If $\mathcal L_u \sim \mathcal L_z$, then each property holding P-a.s. for the
process $z$ must hold also for the process $u$ and vice versa. 
Also the laws of $u(t;x)$ and $z(t;x)$ are equivalent.  In fact,
$P\{u(t;x)\in \Gamma\}=P^*\{z(t;x)\in \Gamma\}$ for every Borel set
$\Gamma \subset H$.
In this way, if we can prove  easily strong Feller property and
irreducibility for the linear equation, 
these properties will be inherited  by the non linear equation.

However, by \eqref{GF:2}
it does not follow that Novikov condition holds.
Anyway, we can approximate the non linearity $F$ in such a way that
Novikov condition holds for the approximate equation and by this we
obtain
$\mathbb E[\rho^T_{u/z}]=1$.
The procedure is standard, but
the results available in the literature 
do not apply here. 
For instance, there are 
similar techniques   
in \cite{mr} (but, even if they deal with 
a stochastic Navier--Stokes  equation, the important 
issue there is the existence of  weak solutions;
Girsanov theorem is proved for other stochastic PDE's) 
or \cite{el} (but, even if they deal with 
a stochastic Kuramoto--Sivashinsky equation,
the Novikov condition 
and Girsanov theorem are analyzed in a finite dimensional context).
We point out that in this paper we prove Girsanov theorem for a 1D stochastic
Kuramoto--Sivashinsky equation and  for a modification of the 2D and 3D
stochastic Navier--Stokes equation. Further,
our results give
regularity of strong solutions of equation \eqref{nlin} (we shall
deal with a variety  of spaces $E\subset H$)
and the equivalence of all its transition
functions so to characterize the asymptotic behaviour by means of Doob
theorem.

\bigskip 

We now state a first result on the absolutely
continuity of the measures.
\begin{proposition}\label{u<z}
Assume \eqref{crescita} holds and that for every $x \in E$ 
there exists a unique strong solution
$z$ of equation \eqref{lin} on the time interval $[0,T]$,
 safisfying \eqref{z:2p}.

Then, given $u(0)=x$ 
there exists a unique weak solution $u$ to equation \eqref{nlin}
on the time interval $[0,T]$ and
the law of the process $u$ 
is absolutely continuous with respect to the law of the process $z$
solving \eqref{lin}, with density  given by \eqref{densita}.
\end{proposition}
\proof
Let us define the 
approximating equation by
\begin{equation}\label{eq:trunc}
\left\{
\begin{array}{l}
 du^N(t)+Lu^N(t)dt+\chi^N_t(u^N) F(u^N(t))dt=Gdw(t)\\
u^N(0)=x
\end{array}
\right.
\end{equation}
where for each $N=1,2,\dots$, the truncation function $\chi^N$ is 
defined as follows:
$$
 \chi^N_t(v)=
\begin{cases}
1 &\text{ if }\int_0^t |G^{-1}F(v(s))|_H^2ds \le N \\
0 &\text{ otherwise } 
\end{cases}
$$
Notice that $\chi^N_{\cdot} (z)$ is a progressively measurable process.
Novikov condition
$$
 \mathbb E \left[\exp
      \left(\tfrac 12 \textstyle\int_0^T 
       |G^{-1}\chi^N_s(z) F(z(s))|_H^2ds\right)
   \right]<\infty
$$
now is trivially satisfied, since by the definition of $\chi^N_t$ we have
$$
 \int_0^T |G^{-1}\chi^N_s(z) F(z(s))|_H^2ds\le N
\qquad P-a.s. .
$$
Hence, for any $N=1,2,\dots$
$$
 \mathbb E[e^{V^{T,N}}]=1 ,
$$
where $V^{T,N}=\int_0^T \chi^N_s(z) \,_H\langle G^{-1}F(z(s)),dw(s)\rangle_H
-\frac 12 \int_0^T \chi^N_s(z) |G^{-1} F(z(s))|_H^2ds$,
and by Girsanov theorem we have that 
$\mathcal L_{u^N}\prec \mathcal L_z$ with the  density
$$
\rho^T_{u^N/z}
=\mathbb E[e^{V^{T,N}}|\sigma^T(z)].
$$

Now we want to prove that $\mathbb E[e^{V^T}]=1$, where the exponent is
$V^T=\int_0^T  \!_H\langle G^{-1}F(z(s)),dw(s)\rangle_H
-\frac 12 \int_0^T |G^{-1} F(z(s))|_H^2ds$. 
\\
We know that $\mathbb E[e^{V^{T,N}}]=1$; moreover
\begin{alignat*}{1}
 \mathbb E[e^{V^{T,N}}]&= \mathbb E[\chi^N_T(z)e^{V^{T,N}} ] 
 +\mathbb E[(1-\chi^N_T(z)) e^{V^{T,N}}]
 \\
 &=\mathbb E[\chi^N_T(z)e^{V^T} ] +P\{\chi^N_T(z)=0\}.
\end{alignat*}
By monotone convergence, 
$\lim_{N\to\infty}\mathbb E[\chi^N_T(z)e^{V^T} ] =\mathbb E[e^{V^T}]$. 
On the other hand, 
$\lim_{N\to\infty} P\{\chi^N_T(z)=0\}=
\lim_{N\to\infty} P\{\int_0^T |G^{-1} F(z(s))|_H^2 ds>N \}=0$.

Therefore $\mathbb E[e^{V^T}]=1$ so that $\mathbb E[e^{V^T}|\sigma^T(z)]$ 
is a probability density. Then, as explained before, 
the probability measure 
$dP^*=\rho^T_{u/z}dP$
(with $\rho^T_{u/z}$ given by \eqref{densita})
defines   a weak solution to equation \eqref{nlin}. 

Uniqueness (in law)  of $u$ is a consequence of uniqueness of $z$ and $\mathcal
L_u\prec \mathcal L_z$. 
\hfill $\Box$

\medskip

Now, besides the previous conditions, let us assume
that also  equation \eqref{nlin}  has a unique strong solution $u$,
enjoying the same property \eqref{z:2p} as $z$.
We obtain a stronger result.

\begin{proposition}\label{uez}
If \eqref{crescita} holds and for any $x\in E$ 
both equations \eqref{nlin} and \eqref{lin}
have a unique strong solution on the time interval $[0,T]$
safisfying \eqref{z:2p},
then the laws $\mathcal L_u$ and $\mathcal L_z$
are equivalent and the densities
are given, respectively, $\rho^T_{u/z}$  by \eqref{densita} and 
$$
 \rho^T_{z/u}=  \mathbb E
   \left[\exp\Big(- \textstyle\int_0^T \!_H\langle G^{-1}F(u(s)),dw(s)\rangle_H
   -\tfrac 12  \textstyle \int_0^T  |G^{-1} F(u(s))|_H^2ds \Big)
   \big|\sigma^T(u) \right].
$$
\end{proposition}
\proof
According to the previous proposition, we know that 
$\mathcal L_u \prec \mathcal L_z$. On the other hand, 
interchanging the r\^ole  of
$u$ and  $z$, again Proposition \ref{u<z}
provides that $\mathcal L_z\prec \mathcal L_u$.
Therefore, they are mutually absolutely continuous, i.e. equivalent.
\hfill $\Box$

\medskip
As a consequence, also the laws of $z(t;x)$ and $u(t;x)$ are equivalent.
Before stating the last result, we need to recall some definitions.
A Markov  process $u$ is said to be strongly Feller 
in $E$ at time $t>0$ if $P_t$ maps $B_b(E)$ into 
$C_b(E)$, where 
$(P_t\phi)(x):=\mathbb E[\phi(u(t;x))]$; and irreducible in $E$ at time $t>0$
if $P(t,x,\Gamma)>0$ for any $x\in E $, $0\neq \Gamma\subset E$ open,
where $P(t,x,\Gamma):=P\{u(t;x)\in \Gamma\}$.

\begin{corollary}\label{il cor}
Under the assumptions of Proposition \ref{uez},
the process $z$ is strongly Feller and irreducible  if and only if 
so is the process $u$.
\end{corollary}

In the next sections, we shall study first the linear equation 
so to check condition \eqref{z:2p} 
and then estimate \eqref{crescita}.

\section{The 1D stochastic Kuramoto--Sivashinsky equation}
\label{1DeqKS}
We refer to \cite{fe} for the abstract setting, in which the 
stochastic Kuramoto--Sivashinsky equation in written as
\begin{equation}\label{sks}
\left\{
\begin{array}{l}
du(t)+ [ \nu A^2 u(t)-Au(t)+B(u(t),u(t))]\ dt=A^\gamma dw(t)\\
u(0)=x
\end{array}
\right.
\end{equation}
and the linear equation associated is
\begin{equation}\label{eqOU}
\left\{
\begin{array}{l} dz(t)+[\nu A^2 z(t) -A z(t)
      +a z(t)]\ dt = A^\gamma dw(t)\\
z(0)=x
\end{array}
\right.
\end{equation}
The unknown $u$ 
can be interpreted as a one-dimensional velocity field
in a compressible fluid (see \cite{uso}).

With respect to the setting of Section \ref{leggi}, we have that 
the linear operator is
$$Lu=\nu A^2u-Au+au$$
with $a>0$   large enough and $\nu >0$,
and the non linear operator  is $$F(u)=B(u,u)-au.$$ 
The operator $G$ in front of the Wiener process is taken of the form 
$A^\gamma$ ($\gamma \in \mathbb R$).
$w$ is a cylindrical 
Wiener process in $H$ on a probability space $(\Omega, F , P)$;
$\{ {F}_t \}_{t \in [0,T]}$ is 
the canonical filtration associated to the Wiener process.

The functional spaces are (given $L>0$, so the spatial domain is
$[-\frac L2,\frac L2]$)
\begin{equation*}
\begin{split}
 &H
 =\{u=u(\xi) \in L^2(-\tfrac L2, \tfrac L2):
          \textstyle\int_{-L/2}^{L/2} u\ d\xi=0\},
\\
 &E
 =D(A^\theta) \text{ for some  } \theta> 0,
\end{split}
\end{equation*}
where
$$
\begin{array}{l}
Au=-u^{\prime\prime}\\
D(A)=H\cap \{u=u(\xi) \in H^2(-\tfrac L2, \tfrac L2):
u(-\tfrac L2)=u(\tfrac L2),
u^\prime(-\tfrac L2)=u^\prime(\tfrac L2)\}.
\end{array}
$$
The operator $A$
is a strictly positive unbounded self-adjoint  operator in $H$, whose
eigenvectors $\{e_j\}_{j=1}^\infty $ form a complete orthonormal basis
of the space $H$.
The powers $A^\theta$ are defined for any $\theta \in \mathbb R$: if
$Ae_j=\lambda_j e_j$ then 
$A^\theta v =\sum_j \lambda_j^\theta \langle v,e_j\rangle e_j$, 
$D(A^\theta)=\{ v=\sum_j v_j e_j:\sum_j \lambda_j^{2\theta} v_j^2<\infty\}$.
Moreover, $\lambda_j \sim j^2$ as $j \to \infty$.

The operator $-(\nu A^2-A+a)$ generates in $H$ (and in any $D(A^\beta)$)
 an analytic semigroup
of negative type of class $C_0$.

The operator $B$ is the  bilinear operator defined by
$$
 B(u,v)=u v^\prime.
$$
For instance, $B$ maps $D(A^{1/2})\times D(A^{1/2})$ into $H$; other
domains
of definition of $B$ are given in \cite{fe}.

First, let us consider the linear equation. 
We are interested in the regularity of the solution $z$ and in the
asymptotic behaviour for $t\to \infty$.
For this, we denote by $R(t,x,\cdot)$ the transitions functions
for \eqref{eqOU}, i.e.
$R(t,x,\Gamma)=P\{z(t;x)\in \Gamma\}$, and by 
$R_t$ the  Markovian semigroup , i.e.
$(R_t\phi)(x)=\mathbb E[\phi(z(t;x))]$.
We say that a measure $m$ is invariant for equation \eqref{eqOU} if
$\int R_t \phi \ dm=\int \phi \ dm$ 
for every $t \ge 0, \phi \in C_b(D(A^\theta))$.
We collect the 
results
in the following proposition.
\begin{proposition}\label{lineareKS}
If $\theta + \gamma <\frac 34 $, then for any $x \in D(A^\theta)$
equation \eqref{eqOU} has a unique strong solution $z$ such that
\begin{equation}\label{ksz:2p}
 \mathbb E \|z\|^{2p}_{C([0,T];D(A^\theta))}<\infty 
\end{equation}
for any $p\ge 1$and $T<\infty$;
this is a Markov process, strongly Feller and irreducible in
$D(A^\theta)$
for any $t>0$. 
The Gaussian measure $\mu_l=
\mathcal N(0,\frac 12 A^{2\gamma}[\nu A^2-A+a]^{-1})$
is the unique invariant measure, 
all transition functions $R(t,x,\cdot)$ are equivalent to
$\mu_l$ 
and
\begin{equation}\label{limiteR}
 \lim_{t \to +\infty} R_t\phi(x)=\int \phi\ d\mu_l, 
\end{equation}
$$
 \lim_{t \to +\infty} R(t,x, \Gamma)= \mu_l(\Gamma)
$$
for any $x \in D(A^\theta),\phi \in C_b(D(A^\theta))$
 and Borel set $\Gamma \subset D(A^\theta)$.
\end{proposition}
\proof
From (3.10) in \cite{fe}, we know that,
given $x \in D(A^\theta)$, if $\theta + \gamma <\frac 34 $
equation \eqref{eqOU} has a unique strong solution
$z$ 
$$
 z(t)=e^{-(\nu A^2-A+a)t}x 
 +\int_0^t e^{-(\nu A^2-A+a)(t-s)} A^\gamma dw(s)
$$
whose paths are, $P$-a.s.,  in $C([0,T];D(A^\theta))$. This is a
Markov process; many of its properties are easy to check, since the
semigroup $\{ e^{-(\nu A^2-A+a)t}\}_{t \ge 0}$
and the covariance of the noise are diagonal operators and 
commute.
 
We recall the basic steps for checking the regularity
of $z$ (the result follows rigorously,
e.g., from \cite{dpz}, Chapter 5, and is proved in \cite{fe}): 

$\; |A^\theta e^{-(\nu A^2-A+a)t}x|_H\le |A^\theta x|_H 
          \qquad \forall t\ge 0$
\begin{equation*}
\begin{split}
 \mathbb E\big|
       \int_0^t A^\theta e^{-(\nu A^2-A+a)(t-s)} A^\gamma dw(s)\big|_H^2
& = \mathbb E\big| \sum_{j=1}^\infty \lambda_j^{\theta+\gamma}
 \int_0^t e^{-(\nu\lambda_j^2-\lambda_j+a)(t-s)} \ d\beta_j(s) \ e_j
 \big|_H^2
\\
& =\sum_{j=1}^\infty \lambda_j^{2(\theta+\gamma)}
 \int_0^t e^{-2(\nu\lambda_j^2-\lambda_j+a)(t-s)}\ ds
\\& \le
 \sum_{j=1}^\infty 
 \frac{\lambda_j^{2(\theta+\gamma)}}{2(\nu\lambda_j^2-\lambda_j+a)}
 \qquad \forall t>0.
\end{split}
\end{equation*}
The last series is convergent if $\theta+\gamma<\frac 34$, since
$\lambda_j \sim j^2$ as $j \to \infty$.
\\
According to Burkholder-Davis-Gundy inequality, the second estimate
provides that inequality \eqref{ksz:2p} holds for any $p$. 

The result on the invariant measure is obtained as in \cite{dpz},
Chapter 11. Actually, the result is trivial if we work first 
on each component $z_j$ and then we recover the infinite dimensional
result for $z$ ($z(t)=\sum_{j=1}^\infty z_j(t) e_j$). 
Indeed, each component $z_j$ 
satisfies
$$
  dz_j(t)+ [\nu \lambda_j^2 -\lambda_j +a]z_j(t)\ dt 
 = \lambda_j^\gamma d\beta_j(t)  , \qquad z_j(0)=x_j;
$$ 
its law is $\mathcal N\big(e^{-(\nu\lambda_j^2-\lambda_j+a)t}x_j,
\frac 12 \frac{\lambda_j^{2\gamma}}{\nu \lambda_j^2-\lambda_j+a}
(1-e^{-2(\nu \lambda_j^2-\lambda_j+a)t})\big)$ 
and for $t \to +\infty$ the density of this Gaussian measure converges
to the density of the Gaussian measure $\mathcal 
N\big(0,\frac 12 \frac{\lambda_j^{2\gamma}}{\nu  \lambda_j^2-\lambda_j+a}
\big)$, which is 
the unique stationary measure.
Therefore,
equation \eqref{eqOU} has a unique invariant measure;
this  is the Gaussian measure with mean 0 and
covariance operator $Q_\infty=\frac 12 A^{2\gamma}[\nu A^2-A+a]^{-1}$;
 
It is easy to check that $\int |A^\theta x|^2_H d\mu_l(x)<\infty$ and
that $\mu_l(\Gamma)>0$ for any open and non empty set 
$\Gamma \subset D(A^\theta)$ .

We expect that irreducibility and strong Feller property hold,  
because the noise acts on all directions $e_j$ and the semigroup
$e^{-(\nu A^2-A+a)t}$
makes  $z$ depending very
regularly on the initial data $x$.

As far as strong Feller property is concerned, by \cite{dpz}
(Chapter 9) we know that condition 
$Ran(Q_t^{1/2})\supset Ran(e^{-(\nu A^2-A+a)t})$
is equivalent to the strong Feller property, where 
$Q_t$ is the covariance operator of the Gaussian random variable 
$z(t;x)$.
Since
$Q_t=\frac 12 A^{2\gamma}[I-e^{-2(\nu A^2-A+a)t}][\nu A^2-A+a]^{-1}$
and for $t>0$ the range of the operator 
$e^{-(\nu A^2-A+a)t}$ is contained in any space $D(A^\beta)$ for $\beta>0$,
we see that this condition  is trivially satisfied.

According to Theorem 11.13 in \cite{dpz}, \eqref{limiteR}
holds and all the transition measures
$R(t,x,\cdot)$ are absolutely continuous with respect to
$\mu_l$. Irreducibility comes straightforward.
Let us point out that in the proof of this theorem, it is shown also that
the law of $z(t;x)$ is equivalent to the law of $z(s;y)$ for
any $t,s>0$ and $x,y\in D(A^\theta)$; 
actually, this follows by Feldman-Hajek theorem,
which is easy to verify in this case of diagonal operators.
\hfill$\Box$

\medskip
To set our problem as in Section \ref{leggi}, we have to fix some space 
$E=D(A^\theta)$.
The interesting spaces are $D(A^\theta)$ for $\theta\ge 0$: $D(A^0)=H$
is the basic space of finite energy  and,  for $\theta>0$,
$D(A^\theta)$ is a subspace
of $H$.
In practise, given $\theta\ge 0$ we choose $\gamma$ as big as possible
($\gamma <\frac 34 -\theta$)
so to make to weakest assumption on the covariance of the noise.
Or, given $\gamma<\frac 34$ (the limitation is
due to $\theta \ge 0$), we choose $\theta$  as big as possible
($\theta<\frac 34 -\gamma$).
Decreasing $\gamma$, the operator $A^\gamma$ is ''more regular'' (in
the sense that, for instance, $A^\gamma$ is a bounded operator for
$\gamma\le 0$) and this stronger assumption provides a more regular
solution $z$ with paths in $C([0,T];D(A^\theta))$.

Now, we deal with estimate \eqref{crescita}.
We have the following result.
\begin{lemma} \label{lem-1}
Let parameters $\gamma $ and $\theta$ be chosen as follows:
\begin{equation}
\begin{split}
\text{for } \tfrac 14 < \gamma < \tfrac 34:\;&\;
\tfrac 38 -\tfrac \gamma 2\le \theta <\tfrac 34 -\gamma
\\
\text{for }0\le \gamma \le \tfrac 14:\;&\;
\tfrac 58 -\gamma\le \theta< \tfrac
34 -\gamma \label{condiz}
\\
\text{for }\gamma <0:\;\;\;\;\;\;\;\;&\;
\tfrac 12 -\gamma
\le \theta <\tfrac 34 -\gamma.
\end{split}
\end{equation}
Then there exists a constant $c$, depending on $\gamma, \theta$ and $a$, such that
$$
 |A^{-\gamma}[B(v,v)-av]|_H\le 
 c \left(1+|A^\theta v|_H^2\right) \qquad \forall v \in D(A^\theta).
$$
\end{lemma}
\proof
Notice that \eqref{condiz} imply  the bounds $\gamma <\frac 34 $, $\theta>0$ and 
$\theta +\gamma<\frac 34$.
The non linear  term is estimated as follows:
\begin{align}
\label{quattro}
&|A^{-\gamma}B(v,\tilde v)|_H \le C_1 |A^{\frac 38 -\frac \gamma 2}v|_H
                                      |A^{\frac 38 -\frac \gamma 2}\tilde v|_H
  &\quad \text{ if }\; \frac 14 < \gamma < \frac 34 
\\\label{cinque}
& |A^{-\gamma}B(v,\tilde v)|_H\le C_2 |A^{\frac 58 - \gamma}v|_H
                                      |A^{\frac 58 - \gamma}\tilde v|_H
  &\quad \text{ if }\; 0\le  \gamma \le \frac 14 
\\\label{sei}
&|A^{-\gamma}B(v,\tilde v)|_H \le C_3 |A^{\frac 12 -\gamma }v|_H
                                      |A^{\frac 12 -\gamma }\tilde v|_H
  &\quad \text{ if } \;\gamma < 0 
\end{align}
The two first inequalities  come from the proof of 
Lemma 2.2 in \cite{gm}.
The latter is proved in  Proposition 2.1 in \cite{fe}.
By the way, recalling that $B(v_1,v_1)-B(v_2,v_2)=
B(v_1-v_2,v_1)+B(v_2,v_1-v_2)$ by bilinearity, 
the above inequalities show that the operator 
$A^{-\gamma}B(v,v)$ is continuous (hence, measurable) 
in the spaces where it is defined.

Notice that if \eqref{condiz} are satisfied, then $\theta>-\gamma$.
Therefore, choosing $\theta$ as in \eqref{condiz} we get
\[
\begin{split}
|A^{-\gamma}[B(v,v)-av]|_H
&\le
 |A^{-\gamma}B(v,v)|_H + a|A^{-\gamma} v|_H
\\
&\le
 C_4 |A^\theta v|_H^2+a C_5  |A^{\theta}v|_H
\\
&\le 
  C_6\left( 1+|A^\theta v|_H^2\right) .
\end{split}
\]
\hfill$\Box$

\begin{remark}\label{casoH}
The case $\theta=0$ is not included.
Indeed,  we have
$$
 |A^{-\gamma}B(v,v)|_H \le c|v|_H^2
$$
for $\gamma>\frac 34$,
because
\begin{equation*}
\begin{split}
 |\langle B(v,v),x\rangle|
 =|\int_{-L/2}^{L/2}\frac 12 (v^2)^\prime x \,d\xi|
 &=\frac 12 |\int_{-L/2}^{L/2} v^2 x^\prime \,d\xi |
 \\&\le \frac 12|v^2|_{L^1} |x^\prime|_{L^\infty}
 \\&\le c |v|_{L^2}^2 |x^\prime|_{D(A^m)} \text{ for } m>\frac 14
 \\&= c |v|_{L^2}^2 |x|_{D(A^{\frac 12+m})} \text{ for } m>\frac 14.
\end{split}
\end{equation*}
       But the condition $\gamma>\frac 34$ is incompatible with
$\theta+\gamma<\frac 34, \theta=0$.
\end{remark}

Now, we consider equation \eqref{sks}.
Let us denote by  $P(t,x,\cdot)$ its transitions functions.

\begin{theorem}
For every $\gamma <\frac 34$ and  choosing $\theta$ as in
\eqref{condiz}, we have the following result.\\
Given $x\in D(A^\theta)$ 
there exist unique strong solutions of equations 
\eqref{sks} and \eqref{eqOU} on any finite time interval $[0,T]$,
with paths in $C([0,T];D(A^\theta))$, $P$-a.s.. We have 
$\mathcal L_u \sim \mathcal L_{z}$, with
the densities 
\begin{multline*}
\rho^T_{u/z}=\mathbb E\left[ e^{V_+^T} \big|\sigma^T(z)\right] \; \text{
  with }V^T_+=\textstyle\int_0^T \!\!\!_H\langle A^{-\gamma}[B(z(s),z(s))-az(s)],
       dw(s)\rangle_H
\\   -\tfrac12 \textstyle\int_0^T
  |A^{-\gamma}[B(z(s),z(s))-az(s)]|_H^2ds 
\end{multline*}
\begin{multline*}
\rho^T_{z/u}=\mathbb E\left[ e^{V^T_-} \big|\sigma^T(u)\right]\; \text{
  with }
V^T_-=-\textstyle\int_0^T \!\!\!_H\langle A^{-\gamma}[B(u(s),u(s))-au(s)],
       dw(s)\rangle_H
 \\  -\tfrac12 \textstyle\int_0^T |A^{-\gamma}[B(u(s),u(s))-au(s)]|_H^2ds 
\end{multline*}
for any $T>0$.\\
Further, $P(t,x,\cdot)\sim \mu_l$ for any $t>0, x\in D(A^\theta)$, 
where $\mu_l=\mathcal N(0,\frac 12 A^{2\gamma}[\nu A^2-A+a]^{-1})$ is
the unique invariant measure for \eqref{eqOU}.
The process $u$ is strongly Feller  and irreducible
in $D(A^\theta)$ at any time $t>0$.
\\
Finally, there exists only one invariant measure $\mu_{KS}$ for
\eqref{sks} which is equivalent to $\mu_l$. 
\\
$\mu_{KS}$ is ergodic, i.e.
$$
 \lim_{T \to +\infty}\frac 1T \int_0^T \phi(u(t;x)) dt =\int \phi \ d\mu_{KS}
$$
$P$-a.s. for every $x \in D(A^\theta), \phi \in L^1(\mu_{KS})$, 
and strongly mixing, i.e.
$$
 \lim_{t \to +\infty} P(t,x,\Gamma)=\mu_{KS}(\Gamma)
$$
for every $x \in D(A^\theta)$ and Borel set $ \Gamma \subset D(A^\theta)$.
\end{theorem}
\proof
If $\theta$ and $\gamma$ are chosen as in \eqref{condiz}, 
from  Proposition \ref{lineareKS} and Lemma \ref{lem-1} we know that
the assumptions of Proposition \ref{u<z} (with $p=2$ and $E=D(A^\theta)$)
are satisfied.
 This implies that 
for $x \in D(A^\theta)$ equation 
\eqref{sks} has a weak solution $u$ living in 
$C([0,T];D(A^\theta))$ and 
$\mathcal L_u \prec \mathcal L_z$; but Theorem 4.3 in \cite{fe} provides
existence and uniqueness of a strong solution $u$ for any $u(0)\in H=D(A^0)$
and $\gamma<\frac 34$.
Thus, we have the regularity result:
 given $x \in D(A^\theta)$ equation \eqref{sks} 
has a  unique strong solution $u$ with paths in $C([0,T];D(A^\theta))$.
By Proposition \ref{uez} we obtain that $\mathcal L_u \sim 
\mathcal L_z$; moreover, $P(t,x,\cdot)\sim R(s,y,\cdot)\sim \mu_l$
and Corollary \ref{il cor} holds. 
We conclude our proof, bearing in mind  Doob
theorem
for uniqueness of invariant measures (see \cite{d-erg}).
The existence of an invariant measure has been proved in \cite{fe}.
\hfill $\Box$

\medskip
Let us notice that, as far as the regularity of solutions is concerned, 
this result improves  that of Proposition 6.5 in
\cite{fe}, since now we can consider any space 
$D(A^\theta)$ with $\theta>0$. 
However, we are not able to prove the  absolute continuity result in
$H=D(A^0)$, as explained in Remark \ref{casoH}, 
even if we know from 
\cite{fe}
that  for any $u(0)\in H$ there exists a unique solution $u$ such that
$u \in C([0,T];H)$ ($P$-a.s.).

Moreover, the results of this section  hold true if the
operator in front of the Wiener process in equation \eqref{sks} is of
the form $LA^\gamma$, 
where $L$ is an isometry in $H$ (e.g., in \cite{fe} we considered
$LA^\gamma w(t)
=\sum_{j=1}^\infty \lambda_j^\gamma  \beta_j(t) (-1)^{j}
e_{j+(-1)^{j+1}}$; 
this includes interesting cases from the physical point of view
as explained in \cite{fe}).

\section{A modified stochastic Navier--Stokes equation}
       \label{navier-stokes}
Since the quadratic term in the Kuramoto--Sivashinsky equation
is similar to that in the Navier--Stokes equation, the only difference
being that the Navier--Stokes equation is set in spaces of divergence free vectors, 
it is appealing to investigate if Girsanov transform holds for the stochastic
        Navier--Stokes equation.
Unfortunately, the answer is negative.
Anyway, let us analyse this problem modifying the linear part.
Our issue is to determine how to modify the  Navier--Stokes equation
to apply our procedure.

Therefore, instead of the stochastic Navier--Stokes equation 
$$
 du(t) +\left[\nu A u(t)+B\big(u(t),u(t)\big)\:\right]\; dt
 =A^\gamma \; dw(t)
$$
(studied, e.g., in \cite{bens}, \cite{vishik}, \cite{Fla}),
we introduce a modification in the linear part; given any $\alpha\ge 1$
we consider 
\begin{equation} \label{sns} 
    \left\{
     \begin{array}{l}
     du(t) +\left[ \nu A^\alpha u(t)+B\big(u(t),u(t)\big)\:\right]\; dt
          =A^\gamma \; dw(t)\\
     u(0)=x
     \end{array} 
    \right.
\end{equation}
This corresponds to replace the Laplacian $-\Delta$ with
$(-\Delta)^\alpha$ in the Navier--Stokes equations in order to seek 
which values of $\alpha$ provide the absolute continuity of $\mathcal
L_u$
with respect to the law of the 
linear equation associated to \eqref{sns}, which is
the modified stochastic Stokes equation:
\begin{equation} \label{ou} 
    \left\{
     \begin{array}{l}
     dz(t) +\nu A^\alpha z(t) dt =A^\gamma \; dw(t)\\
     z(0)=x
     \end{array} 
    \right.
\end{equation}
In this sense, our analysis reminds that of \cite{ms} to investigate 
for which values of $\alpha$ the modified {\it deterministic} Navier--Stokes
equation
$$
\frac{du}{dt}(t) +\nu A^\alpha u(t)+B\big(u(t),u(t)\big)
 =f(t)
$$
is well posed for $d=3$ (we recall that for $d=2$ there is no need of
modification to get existence and uniqueness of a global solution).

With respect to the setting of Section \ref{leggi}, we have that 
the linear operator is
$$Lu=\nu A^\alpha u$$
with $\nu>0$, $\alpha \ge 1$, 
and the non linear operator  is $$F(u)=B(u,u).$$ 
The operator $G$ in front of the Wiener process is taken of the form 
$A^\gamma$ ($\gamma \in \mathbb R$).
$w$ is a cylindrical 
Wiener process in $H$ on a probability space $(\Omega, F , P)$;
$\{ {F}_t \}_{t \in [0,T]}$ is 
the canonical filtration associated to the Wiener process.

The functional setting is defined as usual (see \cite{rosa}).
The symbols $A$ and $B$ will denote different operators from those of Section
\ref{1DeqKS}, but we use the same symbols because of the analogy
between these quantities in
equations
\eqref{sks} and \eqref{sns}.
\\
For $d=2,3$, let
$\dom$ be the $d$-dimensional torus $\mathbb{R}^d/(2\pi \mathbb{Z})^d$,
i.e.  we
consider our problem on the spatial domain  $[0,2\pi]^d$ with  periodic
boundary conditions.

Set
\begin{align*}
  H&=\{ u=\vec u (\vec \xi) \in [L^2(\dom)]^d : 
 \mbox{ div  }u =0, \gamma_n u 
  \mbox{ periodic }, \int_\dom u \ d\vec \xi=0 \}
 \\
  E&= D(A^\theta) \qquad \text{ for some $\theta>0$ }
\end{align*}
where $\gamma_n u$ is the trace of the normal component of $u$ on $\pa \dom$. 

Let $[\dot{H}^m_p(\dom)]^d, m \in \mathbb N \verb+\+ \{0\}$, be
the space of functions of
 $[H^m_{loc}(\mathbb R^d)]^d $, periodic
with period $[0,2\pi]^d$ and with zero average. 
Then the  Stokes operator is defined as
$$
  Au=- \Delta u, \qquad u \in 
          D(A)=[\dot H^2_p(\dom)]^d \cap H. 
$$
$A$ is a strictly positive unbounded self-adjoint operator in $H$, whose
eigenvectors $\{e_j\}_{j=1}^\infty $ form a complete orthonormal basis
of the space $H$.
The powers $A^\alpha$ are defined for any $\alpha \in \mathbb R$. 
The operator $-A$ generates in $H$ 
(and in any $D(A^\beta)$) an analytic semigroup
of negative type $e^{-tA}$
of class $C_0$. Moreover, $A e_j=\lambda_j e_j$ with 
$\lambda_j \sim j^{2/d}$ as $j \to \infty$.

Now, consider the bilinear operator $B$ 
from $D(A^{1/2})\times D(A^{1/2})$ into $D(A^{-1/2})$ defined as
$$
  \langle B(u,v),z\rangle  \:= 
       \int_{\dom} z \cdot [\left(u \cdot \nabla \right) v] \ d\vec\xi 
        \hspace{2cm} \forall \: u,v,z \in D(A^{1/2}).
$$
By the incompressibility condition we have
$$
    \langle B(u,v),v\rangle \:=0, \hspace{2cm} 
    \langle B(u,v),z\rangle \:=-\langle B(u,z),v\rangle. 
$$
Other domains
of definition of $B$ are given below in \eqref{moltipl}.

First, let us consider the linear equation. Similarly to the previous
section, we have 
\begin{proposition}
If
\begin{equation} \label{conv-z}
 \alpha-2(\theta+\gamma)>\frac d2,
\end{equation}
then for any $x \in D(A^\theta)$
equation \eqref{ou} has a unique strong solution $z$ such that
\begin{equation}\label{nsz:2p}
 \mathbb E \|z\|^{2p}_{C([0,T];D(A^\theta))}<\infty 
\end{equation}
for any $p\ge 1$ and $T<\infty$;
this is a Markov process, strongly Feller and irreducible in
$D(A^\theta)$
for any $t>0$. The 
transition functions $\tilde R(t,x,\cdot)$ are equivalent to
$\tilde \mu_l$ for any $t>0, x\in D(A^\theta)$, where $\tilde \mu_l=
\mathcal N(0,\frac 1{2\nu} A^{2\gamma-\alpha})$ 
is the unique invariant measure,
and
\begin{equation}\label{asimpt}
 \lim_{t \to +\infty} \tilde R_t\phi(x)=\int \phi\ d\tilde\mu_l
\end{equation}
$$
 \lim_{t \to +\infty} \tilde R(t,x, \Gamma)= \tilde \mu_l(\Gamma)
$$
for any $x \in D(A^\theta),\phi \in C_b(D(A^\theta))$
and Borel set $\Gamma \subset D(A^\theta)$.
\end{proposition}
\proof
The solution of equation \eqref{ou} is given by
$$
 z(t)=e^{-\nu A^\alpha t}x 
 +\int_0^t e^{-\nu A^\alpha(t-s)} A^\gamma dw(s).
$$
If \eqref{conv-z} holds, then 
there exists a continuous version with values in $D(A^\theta)$.
Indeed,
the basic estimates  are 

$ \qquad\qquad \big|A^\theta e^{-\nu A^\alpha t} x\big|_H
    \le \big|A^\theta x\big|_H 
          \qquad \forall t\ge 0
$
\begin{equation*}
\begin{split} 
 \mathbb E|\int_0^t A^\theta e^{-\nu A^\alpha(t-s)} A^\gamma dw(s)|_H^2
 & =\mathbb E\big| \sum_{j=1}^\infty \lambda_j^{\theta+\gamma}
           \int_0^t e^{-\nu \lambda_j^\alpha(t-s)} \ d\beta_j(s)\;e_j
      \big|_H^2
\\
& =\sum_{j=1}^\infty  \lambda_j^{2(\theta+\gamma)}
 \int_0^t e^{-2\nu\lambda_j^\alpha (t-s)}\ ds
\\& \le
 \sum_{j=1}^\infty 
 \frac{\lambda_j^{2(\theta+\gamma)}}{2\nu\lambda_j^\alpha}
 \qquad \forall t>0.
\end{split}
\end{equation*}
The last series is convergent if 
\eqref{conv-z} is fulfilled, 
since
$\lambda_j \sim j^{2/d}$ as $j \to \infty$.
According to Burkholder-Davis-Gundy inequality, the second estimate
provides that inequality \eqref{nsz:2p} holds for any $p$. 

The unique invariant measure is the Gaussian measure with mean 0 and
covariance operator $\frac 1{2\nu} A^{2\gamma-\alpha}$;
indeed, each component $z_j$ satisfies 
$$
 dz_j(t) +\nu \lambda_j^\alpha z_j(t) dt =\lambda_j^\gamma \; d\beta_j(t)
 ; \qquad     z_j(0)=x_j
$$
and this equation has only one invariant measure which
is the 1-dimensional Gaussian measure 
$\mathcal N(0,\frac 1{2\nu}\lambda_j^{2\gamma-\alpha})$.

\eqref{asimpt} and the equivalence 
$\tilde R(t,x,\cdot)\sim \tilde \mu_l$ can be shown as
in Proposition \ref{lineareKS}.
\hfill$\Box$

\bigskip

Now, we have to choose the space $E$.
Let us consider $\theta\ge 1$.
Why?
Because the easiest estimate for $B(v,v)$ is in the spaces $D(A^m)$
with $m \ge \frac 12$; indeed, for these values the space $D(A^m)$ is a
multiplicative algebra and therefore
\begin{equation}\label{moltipl}
 \big|  A^{m} B(v,\tilde v)\big|_H 
 \le c_m 
 \big| A^m v\big|_H \big|A^{m+\frac 12} \tilde v\big|_H
\end{equation}
(see, e.g., \cite{rosa}). This estimate 
 shows that in these spaces the operator
$A^{m} B(v,v)$ is well defined and continuous.
In particular
\begin{equation}\label{stimaB}
 \big|  A^{\theta-\frac 12} B(v,v)\big|_H 
 \le c^\prime_\theta
 \big| A^{\theta} v\big|^2_H
 \quad \text{ for } \theta \ge 1.
\end{equation}
To check inequality \eqref{crescita} in our context, the latter result
 suggests to 
set 
$$
 -\gamma =\theta-\frac 12.
$$
 In this case, from \eqref{conv-z} we
know that
the process $z$ will have paths in $C([0,T];D(A^\theta))$ if
$$
 \alpha >\frac d2 +1.
$$ 
\begin{remark}
This condition shows that $\alpha=1$ is not allowed.
That is, our procedure does not work for the Navier--Stokes equation;
only taking $\alpha$ sufficiently large we can prove Girsanov theorem and
 the absolute
continuity of the laws.
In particular, for $d=2$ we require $\alpha >2$ and for $d=3$ we
require $\alpha>\frac 52$.
In the same way we can prove this result of absolute continuity for
 the stochastic 1D Burgers equation if $\alpha>\frac 32$.

It is interesting to compare which values of  $\alpha$ provide
that the Navier--Stokes equation is well posed, 
that is it has a unique global solution.
For the {\sl deterministic} equation,
when $d=2$ there is well posedness for $\alpha=1$ 
whereas when  $d=3$ there is well posedness for $\alpha>\frac 54$
(see \cite{ms}).
For the {\sl stochastic} problem, when $d=2$ it is enough to take $\alpha=1$
(see, e.g., \cite{Fla}, \cite{feU}). We guess
that when $d=3$ there is well posedness again for $\alpha > \frac 54$;
this result will be proved in a future work.

\end{remark}

At this point, we prefer to fix a value of $\theta$; indeed,
there are three quantities  involved in the study of  equation \eqref{sns}:
$\alpha, \gamma, \theta$. 
To get not too involved  relations to determine the ''good'' values 
of these parameters, we  reduce the number of parameters
setting  $\theta=1$.
We point out that all the following results can be obtained in the same
way for any $\theta >1$, because of \eqref{moltipl}. 
However, the technicalities are more
involved for $0 \le \theta<1$ (see also Remark \ref{osstheta0} below).

Having set $-\gamma=\theta-\frac 12$,
the choice $\theta=1$ implies $\gamma=-\frac 12$. For these values of
the parameters, we have a pathwise uniqueness result.
This is stronger that uniqueness in law, which would not need to be
proved, as soon as Girsanov transformation holds; indeed, 
if $\mathcal L_u\prec \mathcal L_z$ then uniqueness
of $z$ implies uniqueness in law of $u$.

\begin{proposition}[Pathwise uniqueness]
For $\gamma=-\frac 12$ and $\alpha >\frac d2+1$, given $x \in D(A)$
any two  $C([0,T];D(A))$-valued strong 
solutions of \eqref{sns} coincide $P$-a.s.
\label{pro-unic}
\end{proposition}
\proof
Let $u_1, u_2$ be two
strong  solutions on the probability space $(\Omega, F,\{F_t\}_{t\ge 0},P)$.
Set $U=u_1 - u_2$.
Then $U$ satisfies, $P$-a.s.,
\begin{equation}\label{unico}
 \frac{dU}{dt}(t) +\nu A^\alpha U(t)+B\big(u_1(t),u_1(t)\big)
 -B\big(u_2(t),u_2(t)\big) =0
\end{equation}
with initial data $U(0)=0$. We proceed pathwise.

By bilinearity, 
$B(u_1,u_1) -B(u_2,u_2)= B(u_1,U) +B(U,u_2)$. 
We multiply both sides of \eqref{unico} by $A^2U(t)$; then
(all the norms are in $H$)
\begin{equation*}
\begin{split}
 \frac 12 \frac{d}{dt} |AU(t)|^2+\nu |A^{1+\frac\alpha 2}U(t)|^2
&=-
  \langle
 B\big(u_1(t),U(t)\big)+B\big(U(t),u_2(t)\big),
 A^2 U(t)\rangle
\\
 &=-\langle A^{\frac 12}\big[ B\big(u_1(t),U(t)\big)+
  B\big(U(t),u_2(t)\big)\big], A^{\frac 32}U(t)\rangle.
\end{split}
\end{equation*}
Using \eqref{moltipl},  we have 
$| A^{\frac 12} [B(u_1,U)+B(U,u_2)]|\le c
[|A u_1| +|A u_2|] |AU|$; thus
\begin{equation*}
\begin{split}
 |\langle A^{\frac 12} 
  \big[ B\big(u_1,U\big)+B\big(U,u_2\big)\big],
    A^{\frac 32}U\rangle|
  &\le 
 c \big[|A u_1| +|A u_2|\big] \ |AU|\  |A^{\frac 32}U|
\\ &
\stackrel{(*)}{\le} 
 c \big[|A u_1| +|A u_2|\big] \ |AU|\  |A^{1+\frac \alpha 2}U| 
\\ &
 \le \frac \nu 2 |A^{1+\frac \alpha 2}U|^2
   + c_\nu \big[|A u_1|^2 +|A u_2|^2 \big]|AU|^2.
\end{split}
\end{equation*}
Hence 
$$
 \frac{d}{dt} |AU(t)|^2+\nu |A^{1+\frac\alpha 2}U(t)|^2 
 \le
 C_7\big[|A u_1(t)|^2+ |A u_2(t)|^2\big] |AU(t)|^2.
$$
In particular
$$
  \frac{d}{dt} |AU(t)|^2
 \le
   C_7 \big[|A u_1(t)|^2+ |A u_2(t)|^2\big] |AU(t)|^2.
$$
Since the paths  $ u_1, u_2 \in C([0,T];D(A))$ and $U(0)=0$, 
Gronwall lemma implies that
$$
 |AU(t)|=0 \qquad \forall t \in [0,T],
$$
that is $u_1(t)=u_2(t)$ for all $t \in [0,T]$.
\hfill $\Box$

\begin{remark}
The estimates of the proof remain valid for any $\alpha\ge 1$; 
in fact, inequality $(*)$ holds for $\alpha\ge 1$. Therefore, we could
have stated the proposition assuming only $\alpha \ge 1$.
This depends
strongly on the choice of $\theta$. We point out that for $\theta <1$
uniqueness in $C([0,T];D(A^\theta))$
can be proved along the same lines, but $\alpha$ must be
larger than 1. 

For example, in the case 
 $\underline{\theta=0}$ we estimate the non linearity by
\begin{equation}\label{caso0}
 |A^{-(\frac 12 +\frac d4+\varepsilon)}B(v,\tilde v)|_H \le c |v|_H
  |\tilde v|_H,
\end{equation}
which holds for any $\varepsilon >0$. 
This is proved by means of the embeddings
$D(A^{\frac 12 +\frac d4+\varepsilon})\subset 
[H^{1+\frac d2 +2\varepsilon}(\dom)]^d$ and 
$[H^{1+\frac d2 +2\varepsilon}(\dom)]^d
\subset
[L^\infty(\dom)]^d$, 
that generalize  the estimate
of Remark \ref{casoH} (proved there for $d=1$).
In the proof of pathwise  uniqueness 
(for $\theta=0, \gamma=\frac 12 +\frac d4+\varepsilon$)
we would use
\begin{equation*}
\begin{split}
 |\langle B\big(u_1,U\big)+B\big(U,u_2\big), U\rangle|
&=
 |\langle A^{-(\frac 12 +\frac d4+\varepsilon)} 
  \big[ B\big(u_1,U\big)+B\big(U,u_2\big)\big],
    A^{\frac 12 +\frac d4+\varepsilon}U\rangle|
\\
&  \le 
 c \big[| u_1|_H +|u_2|_H\big]\ |U|_H\ |A^{\frac 12 +\frac d4+\varepsilon}U|_H.
\end{split}
\end{equation*}
If $|A^{\frac 12 +\frac d4+\varepsilon}U|_H\le c
|A^{\frac \alpha 2}U|_H$, that is if
$ \alpha >1+\frac d2$,
we would get that
$$
 \frac{d}{dt} |U(t)|_H^2+\nu |A^{\frac\alpha 2}U(t)|_H^2 
 \le
 C_8 \big[|u_1(t)|_H^2+ |u_2(t)|_H^2\big] |U(t)|_H^2,
$$
so to conclude by Gronwall lemma that $|U(t)|_H=0$ for all $t\in
[0,T]$.

Hence, we can prove pathwise uniqueness in $C([0,T];D(A^0))$ if
$\alpha>1+\frac d2$. On the other hand, chosen 
 $\theta=0$ and $\gamma=\frac 12 +\frac d4+\varepsilon$ 
so to estimate the quadratic term as in  \eqref{caso0},
it follows that inequality \eqref{conv-z} holds for $\alpha>1+d$.

Summing up,
we have checked that to apply our procedure 
for $\theta=0$ we need a stronger assumption on
$\alpha$: $\alpha>1+d$. 
This is the reason for choosing $\theta\ge 1$ so to make the
minimal assumption on $\alpha$.

\label{osstheta0}
\end{remark}

Here is our main result.
\begin{theorem}
For $\gamma=-\frac 12$ and $\alpha >\frac d2+1$, 
given $x \in D(A)$ there exist unique strong solutions of equations 
\eqref{sns} and \eqref{ou} on any finite time interval $[0,T]$; 
the laws $\mathcal L_u$ and 
$\mathcal L_{z}$
are equivalent as measures on the space $C([0,T];D(A))$.

In particular, the densities are 
$$
\rho^T_{u/z}=\mathbb E\left[ e^{+\int_0^T \!_H\langle A^{\frac 12}B(z(s),z(s)),
       dw(s)\rangle_H
    -\frac12 \int_0^T |A^{\frac 12}B(z(s),z(s))|_H^2ds} \big|\sigma^T(z)\right]
$$
$$
 \rho^T_{z/u}=\mathbb E\left[ e^{-\int_0^T \!_H\langle A^{\frac 12}B(u(s),u(s)),
       dw(s)\rangle_H
  -\frac12 \int_0^T |A^{\frac 12}B(u(s),u(s))|_H^2ds} \big|\sigma^T(u)\right]
$$
for any $T>0$.\\
Moreover, $\tilde P(t,x,\cdot)\sim \tilde \mu_l$ for any $t>0, x\in D(A)$, 
where $\tilde \mu_l=\mathcal N(0,\frac 1{2\nu} A^{-1-\alpha})$ 
 is
the unique invariant measure for \eqref{ou}.
In particular, the Markov process $u$ is strongly Feller and
irreducible in $D(A)$ at any time $t>0$; 
hence there exists at most one invariant measure for
\eqref{sns}.
\end{theorem}
\proof For $\theta=1$, $\gamma=-\frac 12$ and $\alpha >\frac d2+1$,
\eqref{conv-z} shows that
the linear equation has a unique strong solution $z$ 
with paths in 
$C([0,T];D(A))$ and satisfying
\eqref{z:2p} for any $p$.
Moreover, by 
\eqref{stimaB} we see  that \eqref{crescita} holds for $p=2$.
According to Proposition \ref{u<z}  we conclude that 
equation \eqref{sns} has a weak solution $u$ living in 
$C([0,T];D(A))$ and satisfying \eqref{z:2p} for any $p$.
This result of weak existence and 
the pathwise uniqueness result of Proposition \ref{pro-unic} 
imply the existence of a strong solution to equation \eqref{sns}
(see, e.g., \cite{ry}, Chapter IX, Th. 1.7). 
By Proposition \ref{uez} we obtain 
$\mathcal L_u \sim \mathcal L_z$ and also $\tilde P(t,x,\cdot) \sim
\tilde \mu_l$.
Then, 
Corollary \ref{il cor} gives strong Feller property and irreducibility
for every $t>0$.
By Doob theorem,  we have uniqueness
of invariant measures for equation \eqref{sns}.

\hfill $\Box$

\begin{remark}
In this section we have assumed periodic boundary conditions 
so to give a meaning 
to terms as $A^{\frac 12}B(z,z)$. The reader can consult
\cite{fe99}  for instance,
to see for which values of $\beta$ the expression $A^\beta B(z,z)$
is well defined when working in a bounded spatial domain 
$\dom\subset \mathbb R^d$,
assuming the velocity vanishes on the boundary $\pa \dom$. 
However, no problem arises
in the periodic case.
\end{remark}

\end{document}